\documentclass[11pt]{article}

\usepackage{latexsym}
\usepackage{amsmath,exscale}
\usepackage{amssymb}
\usepackage{amsfonts}
\usepackage{amsthm}
\usepackage{amscd}
\usepackage{epsfig}
\usepackage{verbatim}
\usepackage{fancybox}
\usepackage{moreverb}
\usepackage{graphicx}
\usepackage{psfrag}
\usepackage[all]{xy}
\usepackage[toc,page]{appendix}

\usepackage[applemac]{inputenc}

\usepackage{color}
\definecolor{marin}{rgb}   {0.,   0.1,   0.9} 
\definecolor{rouge}{rgb}   {0.8,   0.,   0.} 
\definecolor{sepia}{rgb}   {0.8,   0.5,   0.} 
\usepackage[colorlinks,citecolor=marin,linkcolor=rouge,
            bookmarksopen,
            bookmarksnumbered
           ]{hyperref}

\newcommand{\Di}{\displaystyle}
\newcommand{\oo}{\Delta}
\newcommand{\eps}{\varepsilon}	
\newcommand{\feps}{f_{\varepsilon}}

\newcommand{\D}{{\partial}}

\newtheorem{theorem}{Theorem}[section]

\newtheorem{remark}[theorem]{Remark}

\newtheorem{algo}[theorem]{Algorithm}

\def\commutatif{\ar@{}[rd]|{\circlearrowleft}}

\numberwithin{equation}{section}

\title{An exponential integrator for a highly oscillatory Vlasov equation}

\author{Emmanuel Fr\'enod\thanks{Universit\'e Europ\'eenne de Bretagne, LMBA
(UMR CNRS 6205), Universit\'e de Bretagne-Sud \& Inria Nancy-Grand Est, CALVI Project
\& IRMA (UMR CNRS 7501), Universit\'e de Strasbourg, France}
\and
Sever A. Hirstoaga\thanks{Inria Nancy-Grand Est, CALVI Project \& IRMA
(UMR CNRS 7501), Universit\'e de Strasbourg, France}
\and
Eric Sonnendr\"ucker\thanks{Max Planck Institute for Plasma Physics, EURATOM Association,
Boltzmannstr. 2, 85748 Garching, Germany}
}
\date{}

\begin{document}

\maketitle

\begin{abstract}
In the framework of a Particle-In-Cell scheme for some 1D Vlasov-Poisson system
depending on a small parameter, we propose a time-stepping method which is numerically
uniformly accurate when the parameter goes to zero. Based on an exponential time
differencing approach, the scheme is able to use large time steps with respect to the
typical size of the fast oscillations of the solution.
\end{abstract}

\section{Introduction}
\setcounter{equation}{0}
The problem of studying a charged particle beam focused by external electric or magnetic
fields is important in many applications. The motion of such particle beams is governed by
the interactions between the electric field generated by the particles themselves and by
the external focusing electromagnetic field. The modelling framework consists in coupling
a kinetic equation with Maxwell equations. Disregarding the collisions between particles,
the kinetic modelling is performed by means of the Vlasov equation. In this paper we
will consider only non-relativistic long and thin beams. Therefore, instead of studying
the phenomenon by means of
the full Vlasov-Maxwell system, we can use its paraxial approximation (see \cite{Fil-Son06}
for mathematical modelling and numerical simulation of focused particle beams dynamics).
In this framework, the effects of the self-consistent magnetic and electric fields can
be both taken into account by solving a single Poisson equation. Finally, we are led 
to solve a two dimensional phase space Vlasov-Poisson system with a parameter $\eps$.
This small parameter acts on the time variable (in fact the longitudinal variable of a
thin beam, see \cite{Fre-Sal-Son09} for details and physical meaning of the parameter),
producing the highly oscillatory behaviour in phase space. In this framework, the aim of
this paper is to propose a numerical scheme able to study efficiently the evolution of
a beam over a large number of fast periods.

Although more precise alternatives exist, we have chosen in this work to perform the
numerical solution of Vlasov equation by particle methods, which consist in
approximating the distribution function by a finite number of macroparticles.
The trajectories of these particles are computed from the characteristic curves of the
Vlasov equation, whereas the self-consistent electric field is computed on a mesh in the
physical space. This method allows to obtain satisfying results with a small number of
particles (see \cite{BL91}). The contribution of this paper is to propose a new numerical
method for solving the characteristic curves, or equivalently, for computing the
macroparticles' trajectories. Namely, we are faced with solving the following stiff
differential equation
\begin{equation}\label{1steq}
u'(t)=u(t)/\eps + F(t,u(t)),\quad \quad \: u(0)=u_0,
\end{equation}
for several small values of the parameter $\eps$ and where $F$ represents a nonlinear term
which plays the role of the self-consistent electric field. The difficulty arising in the
numerical solution for this equation relies on the ability of the scheme to be uniformly
stable and accurate when $\eps$ goes to zero. Following the survey article
\cite{Hoc-Ost10}, we are encouraged to use efficient numerical methods like exponential
integrators in order to describe the dynamics of equation \eqref{1steq}. The basic idea
behind these methods is to make use of the solution's exact representation given by the
variation-of-constants formula
\begin{equation}\label{voc}
u(t)=e^{t/\eps}\:u_0 +\int_0^te^{(t-\tau)/\eps}F(\tau, u(\tau))\,d\tau.
\end{equation}
Applying this formula from one time step to another has the merit to solve exactly the
linear (stiff) part. Classical numerical schemes fail to capture
the stiff behaviour regardless of the size of the time step with respect to the small
parameter. One may consult \cite{Hoc-Ost10} for construction, mathematical analysis and
implementation of exponential integrators in the two classical types of stiff problems
encapsulated in equation \eqref{1steq}. As a specific implementation, one may cite
\cite{TPLH10}, where, in the context of laser-plasma interactions, an exponential
integrator is used for a Particle-In-Cell (PIC) method in order to model the
high-frequency plasma response.
Once the stiff part is exactly solved, one may use an exponential time differencing (ETD)
method (see \cite{Cox-Matt02}) for the specific numerical treatement of the nonlinear term.
The ETD schemes turn out to outperform many other schemes when treating problems like
\eqref{1steq}; see \cite{Cox-Matt02,Kass-Tref05} for comparisons of ETD against the
Implicit-Explicit method, the Integrating Factor method, etc.
\bigskip 

In the present paper we construct and implement an exponential integrator in order
to solve the characteristics of the highly oscillatory Vlasov-Poisson problem \eqref{VP1d}.
The aim is to use a scheme with large time steps compared to the fast period that arises
from the linear term without loosing the accuracy when the small parameter vanishes. The
novelty of this method is in the numerical approximation of the integral term in
\eqref{voc}. More precisely, when the time step is much larger than the rapid period, the
idea of the algorithm is the following: we first finely solve the ODEs over one fast period
by means of a high-order solver (we have used explicit 4th order Runge-Kutta). Then,
thanks to formula \eqref{voc}, we may compute an approximation of the solution over a large
integral number of periods. We have also found that using a more accurate period, instead
of the period of the solution to the system \eqref{1steq} without the nonlinear term, leads
to more accurate simulations. In addition, we have checked if the scheme gives accurate
solutions starting with an initial condition which lies on the slow manifold or not. We
cite \cite{Boy01} for a ``definition'' of the slow manifold: ``The slow manifold is that
particular solution which varies only on the slow time scale; the general solution to the
ODE contains fast oscillations also.''

The remainder of the paper is organized as follows. Following \cite{Fre-Sal-Son09}, we
briefly recall in Section \ref{sec:axiVP} the paraxial approximation together with the
axisymmetric beam assumption. In Section \ref{sec:PIC-VP} we describe the PIC method for
the Vlasov-Poisson system in which we are interested. Then, Section \ref{sec:eipic} is
devoted to the construction of the new numerical scheme as an exponential integrator for
solving the time-stepping part of the PIC algorithm. Eventually, in Section
\ref{sec:val-numerics}, we implement and test our method on several test cases related
to the Vlasov-Poisson system.

\section{The paraxial model: an axisymmetric Vlasov equation}
\label{sec:axiVP}
The paraxial approximation relies on a scaled Vlasov-Poisson system in a phase space
of dimension four, $2D$ for space variable $\mathbf{x}$ and $2D$ for velocity variable
$\mathbf{v}$. This simplified model of the full Vlasov-Maxwell system is particularly
adapted to the study of long and thin beams of charged particles and it describes their
evolution in the plane transverse to their direction of propagation.

Subject of many research investigations, the paraxial model was derived in a number of
papers, see e.g. \cite{Fil-Son06, Fre-Sal-Son09}. In this work we are interested in
solving numerically a paraxial model with some additional hypotheses, see \eqref{VP1d}
below. The solution of system \eqref{VP1d} is represented by a beam of particles in phase
space. The beam evolves by rotating around the origin in the phase space, and in long
times a bunch forms around the center of the beam from which filaments of particles
are going out. These filaments are difficult to capture with classical numerical methods.

We now introduce the paraxial model that we aim to solve. In the additional axisymmetric
beam assumption (i.e.
invariant beam under rotation in $\mathbb{R}^2$ for $\mathbf{x}$), we are led to change
the $\mathbf{x}$ coordinate in the polar frame. We thus write the model in polar
coordinates $(r,\theta)$, where $r=|\mathbf{x}|$ and $\theta\in[0,2\pi)$ is such that
$x_1=r\cos\theta$ and $x_2=r\sin\theta$.
Then we use new velocity variables $v_r=\mathbf{v}\cdot\mathbf{x}/r$ and
$v_\theta=\mathbf{v}\cdot\mathbf{x}^\perp/r$, where $\mathbf{x}^\perp=(-x_2,x_1)$.
Assuming in addition, as in \cite{Fre-Sal-Son09}, that $\feps$ is concentrated in angular
momentum, i.e. $rv_\theta=0$, the paraxial model becomes
\begin{equation} \label{Vlasov_par_axi}
\left\{
\begin{array}{l}
\cfrac{\D \feps}{\D t} + \cfrac{1}{\eps}\,v_r \cfrac{\D \feps}{\D r} +
\big( \mathbf{E}_r^{\eps} + \Xi_r^{\eps} \big)\cfrac{\D \feps}{\D v_r} = 0, \\ \\
\cfrac1r \cfrac{\D (r\,\mathbf{E}_r^{\eps})}{\D r} = \rho_\eps, \quad 
\rho_\eps(t,r)=\displaystyle\int_{\mathbb{R}}\feps(t,r,v_r)\,dv_r, \\ \\
\feps(t=0,r,v_r) = f_{0}(r,v_r),
\end{array}
\right.
\end{equation}
where the external force $\Xi_r^{\eps}$ writes 
\[
\Xi_r^{\eps}(t,r)=\Big(-\frac{1}{\eps}H_0+H_1\big(\frac{t}{\eps}\big)\Big)\,r,
\]
with $H_0>0$ and $H_1(t)$ some $2\pi$-periodic function with zero mean value.
In this paper we assume that there is no time oscillation in the external field
but only the strong uniform focusing. We thus take $H_0=1$ and $H_1=0$ and the
Vlasov-Poisson system in which we are interested writes
\begin{equation} \label{VP1d}
\left\{
\begin{array}{l}
\cfrac{\D \feps}{\D t} + \cfrac{1}{\eps}\,v_r \cfrac{\D \feps}{\D r} +
\Big( \mathbf{E}_r^{\eps} - \cfrac{r}{\eps} \Big)\cfrac{\D \feps}{\D v_r} = 0, \\ \\
\cfrac1r \cfrac{\D (r\,\mathbf{E}_r^{\eps})}{\D r} =
\displaystyle\int_{\mathbb{R}}\feps(t,r,v_r)\,dv_r, \\ \\
\feps(t=0,r,v_r) = f_{0}(r,v_r).
\end{array}
\right.
\end{equation}
In order to test the numerical method that we propose, we first consider two numerically
simpler test cases where the electric field is not issued from Poisson equation, but it
has analytical forms: in the first case $\mathbf{E}_r^{\eps}(t,r)=-r$ and in the second
one $\mathbf{E}_r^{\eps}(t,r)=-r^3$. Unlike the second case, for the first one, we
can analytically compute the solution to \eqref{VP1d}(a).

\section{A Particle-In-Cell method for Vlasov-Poisson system}
\label{sec:PIC-VP}
We solve the Vlasov-Poisson system \eqref{VP1d} by using a Particle-In-Cell (PIC) method
\cite{BL91}. We thus introduce the following Dirac mass sum approximation of $\feps$
\[
\feps^{N_p}(r,v,t)=\sum_{k=1}^{N_p}w_k\delta(r-R_k(t))\,\delta(v-V_k(t)),
\]
where $N_p$ is the number of macroparticles and $\big(R_k(t),V_k(t)\big)$ is the position
in phase space of macroparticle $k$ moving along a characteristic curve of equation
\eqref{VP1d}(a). Therefore, the problem is to find the positions and velocities
$(R_k^{n+1},V_k^{n+1})$ at time $t_{n+1}$ from their values at time $t_n$, by solving 
\begin{equation}\label{PIC_ODEs}
\left\{
\begin{array}{l}
R'(t)= V(t)/\eps, \\
V'(t)= -R(t)/\eps \,+\,E_{\eps}(t, R(t)),\\
R(t_n)=R_k^n,\,V(t_n)=V_k^n.
\end{array}
\right.
\end{equation}
In this case, the standard PIC algorithm writes as follows: 
(1) deposit particles on a spatial grid, leading to the grid density; 
(2) solve Poisson equation on the grid, leading to the grid electric field;
(3) interpolate the grid electric field in each particle;
(4) push particles with the previously obtained electric field.

The first three steps are classically treated. The first one deals with the computing
of the grid density by convoluting $\rho_{\eps}^{N_p}$ defined by
\[
\rho_\eps^{N_p}(r,t)=\sum_{k=1}^{N_p}w_k\delta(r-R_k(t)),
\]
with a first order spline (this corresponds to the Cloud-in-Cell method in \cite{BL91}). 
Then, we solve Poisson equation \eqref{VP1d}(b) on a uniform one-dimensional grid by using
finite differences. In our case this amounts to only discretize some space integral.
We have done this by using the trapezoidal rule. As for the third step, we used the same
convolution function as for the deposition step in order to get the particle electric
field (this corresponds to a linear interpolation of the grid electric field on each cell).

Eventually, the major issue is the fourth step of the PIC algorithm, consisting in the
numerical integration of system \eqref{PIC_ODEs}. Here is the main focus of this paper,
taking into account that we want to propose a stable and accurate scheme using large time
steps with respect to the fast oscillation. To this end, we introduce in the next section
a method based on exponential time differencing.

\section{The exponential integrator for the Particle-In-Cell method}
\label{sec:eipic}
We now describe the exponential numerical integrator that we have implemented to solve
the fourth step of the PIC algorithm. We first write down the exponential time
differencing (ETD) method in the case of the stiff ODE system we are interested in.
Then, in Section \ref{sec:etd-pic}, we develop an algorithm based on the exponential time
differencing in the framework of the PIC method.

\subsection{Exponential time differencing for a highly oscillatory ODE}
The so-called exponential time differencing scheme arose originally in the field of
computational electrodynamics but has been reinvented many times over the years (see
\cite{Cox-Matt02} and the references therein). We take details of the ideas behind the
various ETD schemes from the comprehensive paper by Cox and Matthews \cite{Cox-Matt02}.

Recall the stiff system of ODEs that we have to solve:
\begin{equation}\label{stiff_ODE}
\left\{
\begin{array}{l}
R'(t)=\;\Di \frac{1}{\eps}\, V(t), \\
V'(t)=-\Di \frac{1}{\eps}\, R(t) \,+\,E_{\eps}(t, R(t)),
\end{array}
\right.
\end{equation}
with some initial condition $(R_0,V_0)$. In this section we assume that the electric field
$E_\eps$ is given. As exposed in \cite{Cox-Matt02}, the stiffness
comes from the two scales on which the solution evolves: the rapid oscillations due to
the linear term and a slower evolution due to the nonlinear (electric) term.
Thus, while any explicit scheme is limited to a small time step, of order $\eps$, a fully
implicit one requires nonlinear problems to be solved and is therefore slow.
A suitable time-stepping scheme for \eqref{stiff_ODE} should be able to avoid the small
time steps when treating the stiffness. The essence of the ETD methods is to solve the
stiff linear part exactly and to derive appropriate approximations when integrating
numerically the slower nonlinear term  \cite{Cox-Matt02,Kass-Tref05}.

To derive the exponential time differencing (ETD) method for this system we first apply
$r(-t/\eps)$ to \eqref{stiff_ODE} and then integrate the obtained equation from $t_n$ to
$t_{n+1}=t_n+dt$ to deduce that
\begin{equation}\label{etd-scheme}
\left(\begin{array}{c}
R_{n+1} \\
V_{n+1}
\end{array}\right)
= r\Big(\frac{dt}{\eps}\Big)\left(\begin{array}{c}
R_n \\
V_n
\end{array}\right)
\;+ \; r\Big(\frac{dt}{\eps}\Big)\;\int_{t_n}^{t_{n+1}}r\Big(\frac{t_n-\tau}{\eps}\Big)
\left(\begin{array}{c}
0 \\
E_{\eps}(\tau, R(\tau))
\end{array}\right) d\tau,
\end{equation}
where 
\begin{equation}
r(\tau) = \left(
\begin{array}{cc}
\cos\,(\tau) & \sin\,(\tau) \\
-\sin\,(\tau) & \cos\,(\tau)
\end{array}
\right).
\end{equation}
This formula is exact. In addition, it is also useful to write \eqref{etd-scheme}
by replacing $(t_n,t_{n+1})$ by any couple $(s,t)$ such that $s<t$. More precisely, we have
\begin{equation}\label{etd-scheme-2}
\left(\begin{array}{c}
R(t) \\
V(t)
\end{array}\right)
= r\Big(\frac{t-s}{\eps}\Big)\left(\begin{array}{c}
R(s) \\
V(s)
\end{array}\right)
\;+ \; r\Big(\frac{t-s}{\eps}\Big)\;\int_s^t r\Big(\frac{s-\tau}{\eps}\Big)
\left(\begin{array}{c}
0 \\
E_{\eps}(\tau, R(\tau))
\end{array}\right) d\tau.
\end{equation}
Now the main question is how to derive approximations to the integral term in
\eqref{etd-scheme}. All the ETD schemes are results of this process.
In this spirit, it was shown in \cite{Cox-Matt02}
that ETD methods can be extended to any order by using multistep or Runge-Kutta
methods and explicit formulae for such arbitrary order ETD methods were derived.
In particular, explicit coefficients for ETD Runge-Kutta methods of order up to four
have been computed. The authors also illustrated on several ODEs and PDEs that ETD is
superior over the Integrating Factor and Implicit-Explicit methods, two other
classical schemes able to avoid the small time step.
Nevertheless, in the form written in \cite{Cox-Matt02}, a high-order ETD scheme
(e.g. ETDRK4) suffers from numerical instability as explained in \cite{Kass-Tref05}.
The problem has been solved in \cite{Kass-Tref05} by using a contour integral method for
evaluating the coefficients. Then, this modified ETDRK4 scheme has been tested in
\cite{Kass-Tref05} against five other 4th order schemes on several PDEs, and ETDRK4 has
been found the best in terms of errors. These results encouraged us to use an ETD scheme
in order to solve stiff ODEs.

In this paper we do not use a high-order ETD method as described in the references
above but merely the formula \eqref{etd-scheme} for justifying the derivation of
Algorithm \ref{algo-time}. More precisely,
we adopt a different approach for approximating the integral term in \eqref{etd-scheme},
which is justified by the remark in the next paragraph and by our aim to use very large
time steps compared to the fast oscillations, say $dt=100\,\eps$ when $\eps=10^{-2}$.
This is the core of the method described in the following section.

\begin{figure}[h]
\begin{center}\hspace*{-5mm}
\begin{tabular}{cc}
\includegraphics[scale=0.65]{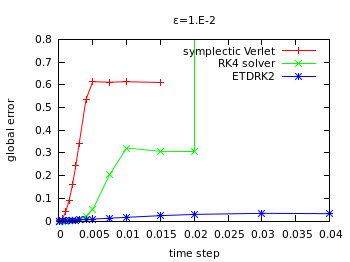} &
\includegraphics[scale=0.65]{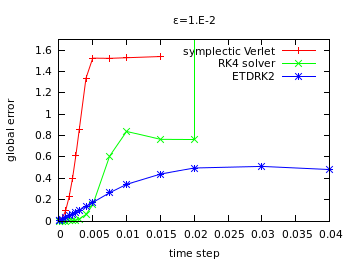}
\end{tabular}
\end{center}
\caption{A nonlinear case: the global Euclidean error at time 3.14 for two initial
conditions: on (left) and off (right) the slow manifold}\label{etdrk2vsrk4sv}
\end{figure}
However, we have done tests with ETDRK2 (see formula $(22)$ in \cite{Cox-Matt02}) against
symplectic Verlet and RK4 schemes, in the case $E_\eps(t,R)=-R^3$ with $\eps=10^{-2}$ (see 
the test case in Section \ref{sec:1nonlinear}). Thus, we have noticed that this ETD
scheme has smaller errors for some range for the
time step and that, unlike its competitors, it is stable when $dt\geq 2\eps$ (see
Fig.~\ref{etdrk2vsrk4sv}). But these simulations also show that when $dt$ is too large
with respect to $\eps$ (e.g. $dt\geq 3\eps$) the global error is significant for some
initial condition (particles off the slow manifold, see Section \ref{sec:val-numerics}).
The reason is that for $dt$ large enough with respect to $\eps$, the nonlinear force term
cannot be accurately taken into account by any high-order Runge-Kutta solver.

\subsection{The ETD PIC method with large time steps}
\label{sec:etd-pic}
In this section we describe our time-stepping scheme starting from the equation
\eqref{etd-scheme}.
The main idea is to remark that the time for one fast grand tour is approximated by
$2\pi\eps$, the period of the solution
to \eqref{stiff_ODE} without nonlinear term. Therefore, since we want to build a scheme
with a time step $dt$ much
larger than the fast oscillation, we first need to find the unique positive integer $N$
and the unique real $\oo\in[0,2\pi\eps)$ such that
\begin{equation}\label{deco-dt}
dt=N\cdot(2\pi\eps)\,+\,\oo.
\end{equation}
Thus the integral term in \eqref{etd-scheme} becomes
\begin{equation}\label{deco:integ}
\int_{t_n}^{t_{n+1}}d\tau=\sum_{j=0}^{N-1}\int_{t_n+2\pi\eps\,j}^{t_n+2\pi\eps\,(j+1)}d\tau+
 \int_{t_{n+1}-\oo}^{t_{n+1}}d\tau,
\end{equation}
that we approximate by
\begin{equation}\label{approx:1}
N\,\int_{t_n}^{t_n+2\pi\eps} d\tau + \int_{t_{n+1}-\oo}^{t_{n+1}}
d\tau.
\end{equation}
In this approximation, we make some assumptions that we develop in Remark \ref{rem:1}.
Computing the integrals 
\begin{equation*}
{\cal I}_1 = \int_{t_n}^{t_n+2\pi\eps} r\Big(\frac{t_n-\tau}{\eps}\Big)
\left(\begin{array}{c}
0 \\
E_{\eps}(\tau, R(\tau))
\end{array}\right) d\tau\;
\textrm{and}\;
{\cal I}_2 = \int_{t_{n+1}-\oo}^{t_{n+1}} r\Big(\frac{t_n-\tau}{\eps}\Big)
\left(\begin{array}{c}
0 \\
E_{\eps}(\tau, R(\tau))
\end{array}\right) d\tau
\end{equation*}
knowing $(R_n,V_n)$ will thus lead to $(R_{n+1},V_{n+1})$. Next we describe this method in
3 steps. The electric field is discretized explicitely, i.e. it can be computed at any
time solving the Poisson equation \eqref{VP1d}(b) with a right-hand side obtained from
depositing the particles on the grid.
\medskip

{\bf 1st step}: Using \eqref{etd-scheme-2} with $s=t_n$ and $t=t_n+2\pi\eps$ and since
$r(2\pi)=\textrm{Id}$, we have 
\begin{equation}\label{approx:I1}
{\cal I}_1 = 
\left(\begin{array}{c}
R(t_n+2\pi\eps) - R(t_n)\\
V(t_n+2\pi\eps) - V(t_n)
\end{array}\right).
\end{equation}
Now, we finely solve \eqref{stiff_ODE} with initial conditions $(R_n,V_n)$ by a 4th order
Runge-Kutta scheme, in order to get an accurate approximation of $\big(R(t_n+2\pi\eps),
V(t_n+2\pi\eps)\big)$. This is done by following the four steps of the standard PIC
algorithm as exposed in Section \ref{sec:PIC-VP}. The time step used in the fourth step
for the RK4 solver is $\sqrt\eps\,\eps$ (we explain this choice in
Section~\ref{sec:val-numerics}).
\medskip

{\bf 2nd step}: We have to calculate an approximation of $\big( R(t_n+N\cdot2\pi\eps),
V(t_n+N\cdot2\pi\eps) \big)$ since it will be usefull in the computations to do in
the 3rd step. Using \eqref{etd-scheme-2} with $s=t_n$ and
$t=t_n+N\cdot2\pi\eps$ we obtain, since $r(N\cdot2\pi)=\textrm{Id}$, 
\begin{equation}\label{interm:approx}
\left(\begin{array}{c}
R(t_n+N\cdot2\pi\eps) \\
V(t_n+N\cdot2\pi\eps)
\end{array}\right)
\approx \left(\begin{array}{c}
R_n\\
V_n
\end{array}\right)\,+\, N\cdot {\cal I}_1,
\end{equation}
where the right integral (the first one in the right-hand side of \eqref{deco:integ})
was approximated as done in \eqref{approx:1} by $N\cdot {\cal I}_1$.
\medskip

{\bf 3rd step}: Now, we compute  ${\cal I}_2$ by using \eqref{etd-scheme-2} with
$s=t_n+N\cdot2\pi\eps=t_{n+1}-\oo$ and $t=t_{n+1}$
\begin{equation}\label{approx:I2}
{\cal I}_2 = r\Big(-\frac{\oo}{\eps}\Big)
\left(\begin{array}{c}
\widetilde R(t_{n+1})\\
\widetilde V(t_{n+1})
\end{array}\right)\; - \;
\left(\begin{array}{c}
R(t_n+N\cdot2\pi\eps) \\
V(t_n+N\cdot2\pi\eps)
\end{array}\right),
\end{equation}
where $\big(\widetilde R(t_{n+1}),\widetilde V(t_{n+1})\big)$ may be found as done above
for the approximation of $\big(R(t_n+2\pi\eps),$\\$V(t_n+2\pi\eps)\big)$: we follow the
steps of the standard PIC algorithm, using for the particles' push a 4th order
Runge-Kutta solver with initial conditions $(R(t_n+N\cdot2\pi\eps),V(t_n+N\cdot2\pi\eps))$
and with total time $\oo$ (as for the 1st step above, we choose a small time step,
$\sqrt\eps\,\eps$).
\medskip

Finally, we replace \eqref{interm:approx} in  \eqref{approx:I2} which we put in
\eqref{etd-scheme}: the term $N\cdot{\cal I}_1$ will cancel and at the end we have 
$$R_{n+1}=\widetilde R(t_{n+1})\quad\text{and}\quad V_{n+1}=\widetilde V(t_{n+1}),$$
meaning that the vector $\big(\widetilde R(t_{n+1}),\widetilde V(t_{n+1})\big)$ calculated
within the 3rd step above is an approximation of the solution at time $t_{n+1}$.
\bigskip

In  summary, writing $dt$ as in \eqref{deco-dt}, the implementation of our
time-stepping algorithm follows three steps:
\begin{algo}\label{algo-time}
Assume that $(R_n,V_n)$ the solution of \eqref{stiff_ODE} at time $t_n$ is given. Then
\begin{enumerate}
\item compute $(R,V)$ at time $t_n+2\pi\eps$ by using a fine Runge-Kutta solver with
initial condition $(R_n,V_n)$.
\item compute $(R,V)$ at time $t_n+N\cdot2\pi\eps$ by the following rule
\begin{equation}\label{interm:approx-2}
\left(\begin{array}{c}
R(t_n+N\cdot2\pi\eps) \\
V(t_n+N\cdot2\pi\eps)
\end{array}\right)
= \left(\begin{array}{c}
R_n\\
V_n
\end{array}\right)\,+\, N
\left(\begin{array}{c}
R(t_n+2\pi\eps) - R_n\\
V(t_n+2\pi\eps) - V_n
\end{array}\right).
\end{equation}
\item compute $(R,V)$ at time $t_{n+1}$ by using a fine Runge-Kutta solver with
initial condition $(R,V)$ obtained at the previous step.
\end{enumerate}
\end{algo}
In the framework of the PIC method, Algorithm~\ref{algo-time} is applied to each particle.
\medskip

We will call the \emph{modified ETD-PIC scheme}, the Algorithm~\ref{algo-time} where
$2\pi\eps$ is replaced in the first two steps with a more accurate fast time for particles.
The reason for this choice is explained in Sections \ref{sec:linear} and
\ref{sec:comm}.

\begin{remark}\label{rem:1}
In the approximation \eqref{approx:1}, we have made the assumptions that the solution's
period does not change significantly in time and the same for the electric field
$E_{\eps}$. Nevertheless, although in this section $E_\eps$ in supposed to be given, in
the case of the Vlasov-Poisson coupling \eqref{VP1d}, $E_\eps$ is the self-consistent
electric field. Therefore, we expect an almost periodic trajectory of the charged
particles to involve a similar
time behaviour for the electric field that they generate. Thus, in such a case, we make
only the assumption that the particles' period does not vary significantly in time.
\end{remark}

\section{Validation of the numerical method}
\label{sec:val-numerics}
We now validate our algorithms in the test cases presented in Section \ref{sec:axiVP}.
We illustrate their numerical performance and desired properties exposed in the
Introduction with global error curves. Thus, we have checked if the method is numerically
stable and accurate when $\eps$ becomes smaller. Since we solve a stiff system, we compute 
the errors for both types of initial particles, on and off the slow manifold (see 
\cite{Boy01,Cox-Matt02}).
In fact, for numerical reasons, we should rather use the designations ``close or far from
the slow manifold'' than ``on or off the slow manifold''.
Nevertheless, in order to be in line with the literature we cite, we keep in this paper the 
designations ``on or off the slow manifold''.

In order to solve the Vlasov-Poisson model \eqref{VP1d}, we choose as initial distribution
function the following semi-Gaussian beam (see \cite{Fre-Sal-Son09})
\begin{equation}
f_0(r,v)= \frac{n_0}{\sqrt{2\pi}\,v_{\rm{th}}}\exp\Big(-\frac{v^2}{2v_{\rm{th}}^2}\Big)
\chi_{[-0.75,0.75]}(r),
\end{equation}
where the thermal velocity is $v_{\rm{th}}=0.0727518214392$ and $\chi_{[-0.75,0.75]}(r)=1$
if $r\in [-0.75,0.75]$ and $0$ otherwise. We implement this distribution using a particle
approximation, with $N_p=10000$ macroparticles with equal weights $w_k=1/N_p$.

\subsection{A linear case}
\label{sec:linear}
As a first test of validation of our scheme we consider a case
where the solution is analytically known. It is the linear case where the self-consistent
electric field in \eqref{stiff_ODE} is given by $E_{\eps}(t, R)= -R$. The solution
satisfying the initial condition $R(s)=R_s$, $V(s)=V_s$ is
\begin{equation}\label{stiff_ODE_1_sol}
\left\{
\begin{array}{l}
R(t)= V_s\Di\frac{1}{\sqrt{1+\eps}}\,\sin\Big(\Di\frac{\sqrt{1+\eps}}{\eps}\,(t-s)\Big)
+R_s\,\cos\Big(\Di\frac{\sqrt{1+\eps}}{\eps}\,(t-s)\Big)\\
V(t)= V_s\,\cos\Big(\Di\frac{\sqrt{1+\eps}}{\eps}\,(t-s)\Big)
-R_s\,\sqrt{1+\eps}\,\sin\Big(\Di\frac{\sqrt{1+\eps}}{\eps}\,(t-s)\Big)
\end{array}
\right.
\end{equation}
for all $t\geq s$. We note that the trajectory in the phase space is an ellipse since we
have $$\forall t\geq s \qquad R(t)^2+\Big(\frac{V(t)}{\sqrt{1+\eps}}\Big)^2=R_s^2+
\Big(\frac{V_s}{\sqrt{1+\eps}}\Big)^2.$$
In this particular case, we can compute exactly the time for one rapid grand tour.
It is \begin{equation} t_p=\frac{2\pi\eps}{\sqrt{1+\eps}} \end{equation}
and this value (and not $2\pi\eps$) was
used in our simulation for the push of particles within Algorithm \ref{algo-time}. In
fact, using $2\pi\eps$ instead of the right rapid time makes our method unstable when the
final time is large enough. The reason is that $2\pi\eps-t_p$ is big enough so that the
accumulated in time errors issued from the rule \eqref{interm:approx-2} lead the particles
to drift outward in the phase space. We also notice that, unlike general case, the rapid
time does not depend on the initial condition. Thus, the beam of particles is rotating in
the phase space without spiraling. This was checked numerically for our scheme.

In addition, in this case, the slow manifold can easily be determined. Knowing that it is
that particular solution which varies only on the slow time scale (\cite{Boy01}), we deduce
that the fast oscillations in the solution given by \eqref{stiff_ODE_1_sol} can be removed
when $R_0=0$ and $V_0=0$. The stationary solution $(R,V)=(0,0)$ is the slow manifold in
this particular case. Numerically, on and off the slow manifold particles were arbitrarly
chosen with $R_0=0.306825$ and $V_0\sim 7\cdot10^{-6}$ (close to $(0,0)$) and with
$R_0=0.748725$ and $V_0\sim 0.142892$, respectively. 

\begin{figure}[h]
\begin{center}\hspace*{-5mm}
\begin{tabular}{cc}
\includegraphics[scale=0.65]{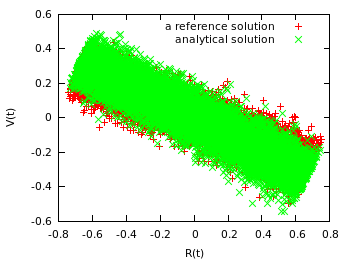} &
\includegraphics[scale=0.65]{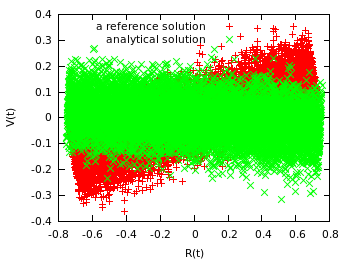}
\end{tabular}
\end{center}
\caption{Computing a reference solution in phase space: use of an explicit Runge-Kutta
solver with $dt=\eps/10$ when $\eps=10^{-5}$ at final time 1.0 (left) and 3.5 (right)}
\label{solRK4bigdt}
\end{figure}

Eventually, this simple case allows us to remark the following. When computing a reference
solution, we push the particles with explicit 4th order Runge-Kutta solver using a
sufficiently small time step. As a result of several numerical experiments, we found the
quite optimal $dt=\sqrt{\eps}\,\eps$. More precisely, we have first noticed the need for
the choice of a uniformly small ratio $dt/\eps$ with respect to $\eps$ in order to get an
accurate solution. Indeed, we have computed a reference solution with a time step
$dt=\eps/10$, when $\eps=10^{-5}$. The outcome of this simulation at final times $t=1.0$
and $t=3.5$ is illustrated in Fig.~\ref{solRK4bigdt}, showing that the time step
$dt=\eps/10$ is not sufficiently small so that good accuracy for the reference solution be
reached. Second, we have experimented different uniformly small time steps, $dt=\eps^{5/4}$,
$\eps^{3/2}$, $\eps^{7/4}$, $\eps^2$. Thus, computing solutions with time steps
$dt\leq\eps^{7/4}$ do not lead to considerably smaller error, for small $\eps$, and in
addition it requires large CPU time. We have finally declared acceptable, at worst of
order $10^{-5}$, the errors obtained with $dt=\eps^{3/2}$, and therefore, this time
step was chosen in the following test cases when computing the reference solution.

\subsection{A nonlinear case}
\label{sec:1nonlinear}
We now take the case of $E_{\eps}(t, R)= -R^3$.
The solution to system \eqref{stiff_ODE} can be writen in terms of Jacobi elliptic
functions and it cannot be put in an analytical form. We therefore compute a reference
solution using a fine explicit Runge-Kutta solver. Unlike the nonlinear case in Section
\ref{sec:vp}, we do not have numerical errors due to the computation of the electric field.
\begin{figure}[h]
\begin{center}\hspace*{-5mm}
\begin{tabular}{cc}
\includegraphics[scale=0.65]{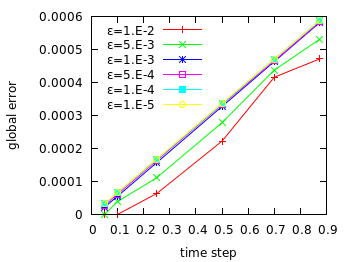} &
\includegraphics[scale=0.65]{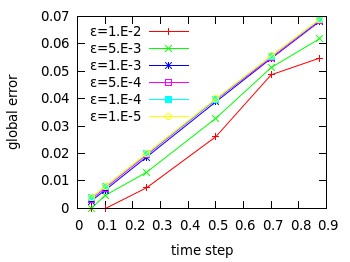}
\end{tabular}
\end{center}
\caption{First nonlinear case: the global Euclidean error at time 3.5 for an initial
particle on (left) and off (right) the slow manifold}\label{R3}
\end{figure}

We have used the initial conditions, on and off the slow manifold, mentioned in the
previous section, since $(R,V)=(0,0)$ is still a stationary solution to problem
\eqref{stiff_ODE}. We have applied Algorithm~\ref{algo-time} for each particle. 
Denoting by $(R(t),V(t))$ the result of the ETD-PIC method and by $(R_{\text{ref}}(t),
V_{\text{ref}}(t))$ the reference solution of problem \eqref{stiff_ODE} for some initial
condition, the global errors (the max value in time of the Euclidean norm of
$(R(t)-R_{\text{ref}}(t),V(t)-V_{\text{ref}}(t))$)
at final time $3.5$ are shown in Fig.~\ref{R3} for several values of $\eps$. The error
curves in the left panel are derived with the above initial condition on the slow manifold.
\begin{figure}[h]
\begin{center}\hspace*{-5mm}
\begin{tabular}{cc}
\includegraphics[scale=0.65]{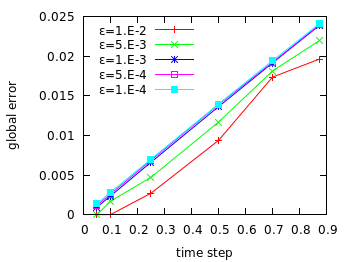} &
\includegraphics[scale=0.65]{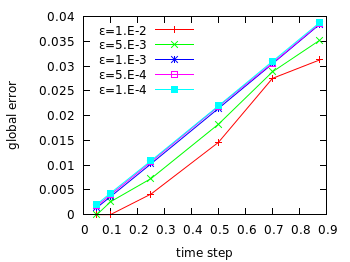}
\end{tabular}
\end{center}
\caption{Vlasov-Poisson case: the global Euclidean error at time 3.5 for an initial
particle on (left) and off (right) the slow manifold}\label{cwithPoiss}
\end{figure}
\medskip

\subsection{The Vlasov-Poisson case}
\label{sec:vp}
No analytical solution to system \eqref{stiff_ODE} is available in the case where
$E_{\eps}$ is the solution of Poisson equation \eqref{VP1d}(b). We solve numerically this
equation by using a trapezoidal formula for the integral in $r$ with $128$ cells for the
space interval. As above, the time step of the RK4 solver for computing the reference
solution is $dt=\sqrt{\eps}\,\eps$ and we have used the same initial conditions, on and
off the slow manifold. The global errors of the ETD-PIC method are shown in
Fig.~\ref{cwithPoiss} for different values of $\eps$.

\begin{remark}
\begin{enumerate}
\item From Figs.~\ref{R3} and \ref{cwithPoiss} one can see the announced property of
the scheme, that is the uniform accuracy when $\eps$ goes to $0$. We also observe
that the order of accuracy decreases \text{uniformly} when $\eps$ goes to $0$.

\item Note the very large time steps with respect to $\eps$ that the ETD scheme allows
to use. For instance, when the time step is $0.7$ this means that it goes from $70\,\eps$
when $\eps=10^{-2}$ to $70000\,\eps$ when $\eps=10^{-5}$.
\end{enumerate}
\end{remark}
\begin{figure}[h]
\begin{center}\hspace*{-5mm}
\begin{tabular}{cc}
\includegraphics[scale=0.65]{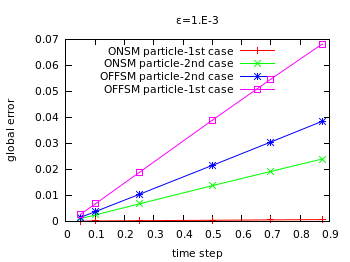} &
\includegraphics[scale=0.65]{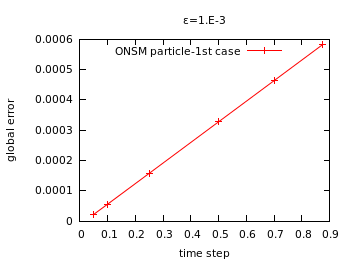}
\end{tabular}
\end{center}
\caption{Global errors at time 3.5 for initial particles
on and off the slow manifold in both nonlinear cases (1st case is the one in
Section \ref{sec:1nonlinear}, 2nd case is the Vlasov-Poisson case)}\label{rototo}
\end{figure}
\subsection{Comments}
\label{sec:comm}
~~~  First, we notice that for each nonlinear case above, the errors are larger when
the initial condition is off the slow manifold.
Then,  when $\eps=10^{-3}$, we can see in Fig.~\ref{rototo} that the difference between
errors for off and on the slow manifold particles is more significant in the first
nonlinear test case.

Second, still when $\eps=10^{-3}$, we have compared the two nonlinear test cases.
When simulation is done with an initial particle on the slow manifold (ONSM), the
error is much more significant for the second nonlinear case, as expected because of
the additional numerical errors due to Poisson solver (see Fig.~\ref{rototo}).
Surprisingly, for an initial particle off the slow manifold (OFFSM) the error behaves
conversely, although they are of the same order of magnitude. 

We think the reason for both comments above is in the use of $2\pi\eps$ as the fast
period for all the particles, which is more or less accurate. Precisely, we have
experimentally found that while the period of an ONSM
particle in the first nonlinear case is very close to $2\pi\eps$ (thus justifying
the much smaller error at the right of Fig.~\ref{rototo}), in the second nonlinear
case this period is rather close to $2\pi\eps + 3\eps^2$. In addition, for an OFFSM
particle our experiments lead to fast times close to $2\pi\eps-\eps^2$ in the first
nonlinear case and to $2\pi\eps + \eps^2$ in the second nonlinear case. Even if these
numbers seem to be small to induce such significant errors, we recall that in the linear
test case above, $2\pi\eps-t_p$ is very close to $3\eps^2$ and the use of a not
accurate fast time for particles leads to an unstable simulation. That is why we test in
the following section {\em modified ETD-PIC schemes}, where we use a more accurate fast
time than $2\pi\eps$. 

\begin{figure}[h]
\begin{center}\hspace*{-5mm}
\begin{tabular}{cc}
\includegraphics[scale=0.65]{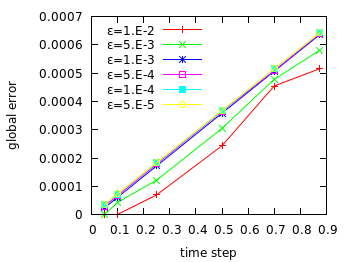} &
\includegraphics[scale=0.65]{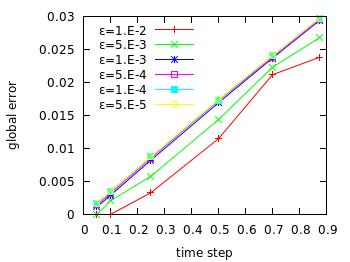}
\end{tabular}
\end{center}
\caption{First nonlinear case, Algorithm \ref{algo-time} with the particles mean period:
the global Euclidean error at time 3.5 for an initial particle on (left) and off (right)
the slow manifold}\label{R3withmeanper}
\end{figure}
\begin{figure}[h]
\begin{center}\hspace*{-5mm}
\begin{tabular}{cc}
\includegraphics[scale=0.65]{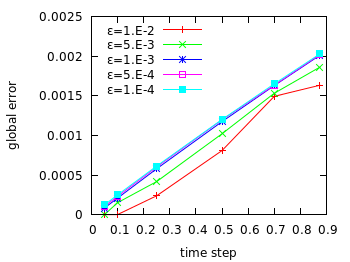} &
\includegraphics[scale=0.65]{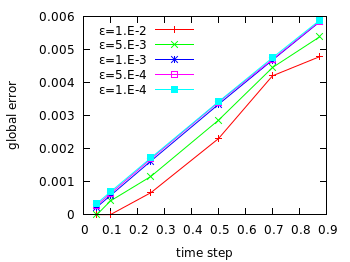}
\end{tabular}
\end{center}
\caption{Vlasov-Poisson case, Algorithm \ref{algo-time} with the particles mean period:
the global Euclidean error at time 3.5 for an initial particle on (left) and off (right)
the slow manifold}\label{cwPoisswmeanper}
\end{figure}

\subsection{Validation of the modified ETD-PIC numerical scheme}
Finally, we present the numerical results obtained with Algorithm \ref{algo-time} where
$2\pi\eps$ was replaced with a more accurate period. More precisely, we first compute
numerically an approximation of the period for each particle, based on the fact that
particles trajectories in time are sinusoid-like almost periodic functions. The
trajectories are computed with the same Runge-Kutta solver used for the reference
solution's calculation. Thus, the particles are pushed during some very small time until
all the particles reach their trajectory's third extremum. The criterion for finding these
extrema is the velocity's change of sign. We then stated that each particle's period is
the time interval between the first and the third extremum. Eventually, we compute the
mean of these periods. It is by this mean that we replace $2\pi\eps$ in the first two
steps of Algorithm \ref{algo-time}. Summarising, the implementation of the modified
ETD-PIC scheme starts with the finding of the mean period and then uses this value in
Algorithm \ref{algo-time} all along the simulation.

Our experiments show that, at best, in the coupling with Poisson case, the error is smaller
by a factor of $10$ than that computed with period $2\pi\eps$ (see Figs. \ref{cwithPoiss}
and \ref{cwPoisswmeanper}). For the nonlinear case in Section \ref{sec:1nonlinear}, as
expected (see Section \ref{sec:comm}), the errors are slightly larger than those computed
with period $2\pi\eps$ for an initial ONSM particle, but smaller for an initial OFFSM
particle.

\begin{figure}[ht]
\begin{center}\hspace*{-5mm}
\begin{tabular}{cc}
\includegraphics[scale=0.65]{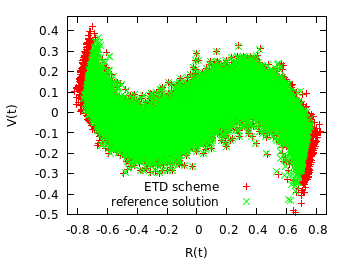}%2pieps_R3.jpg en 0.25}&% the file with 2\pi\eps for period
\includegraphics[scale=0.65]{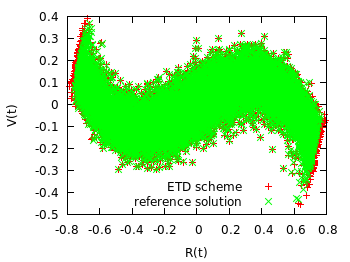}%jpg  en 0.25}
\end{tabular}
\end{center}
\caption{First nonlinear case, $\eps=10^{-4}$: phase space solution computed with a time step
$8750\:\eps$ at final time 3.5, using for particles period $2\pi\eps$ (at left) and the mean
period (at right)}\label{compar_sol_meanper_R3}
\end{figure}

\begin{figure}[ht]
\begin{center}\hspace*{-5mm}
\begin{tabular}{cc}
\includegraphics[scale=0.65]{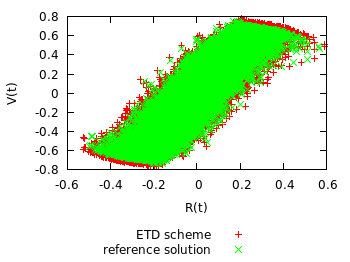}%_2pieps.jpg} &  en 0.25
\includegraphics[scale=0.65]{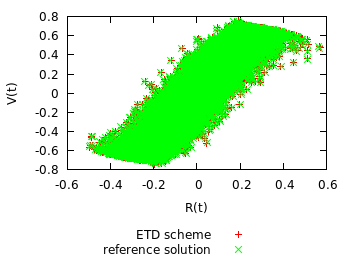}%.jpg}  en 0.25
\end{tabular}
\end{center}
\caption{Vlasov-Poisson case, $\eps=10^{-4}$: phase space solution computed with a time step
$8750\:\eps$ at final time 3.5, using for particles period $2\pi\eps$ (at left) and the mean
period (at right)}\label{compar_sol_meanper_cwP}
\end{figure}

At the end, we also illustrate the results of ETD-PIC scheme vs. {\em modified ETD-PIC
scheme} by representing the solution to Vlasov equation with an electric field, first
given as in Section \ref{sec:1nonlinear} (Fig. \ref{compar_sol_meanper_R3}) and second as
solution of the Poisson equation (Fig. \ref{compar_sol_meanper_cwP}). These beams of particles
were obtained with the larger time step used for the errors calculus, $dt=0.875$, meaning
only $4$ iterations for the ETD-PIC schemes to reach $t=3.5$.

\section{Conclusion}
In this paper we have introduced a new numerical scheme for solving some stiff (highly
oscillatory) differential equation. This scheme is based on exponential time differencing
and can accurately handle large time steps with respect to the fast oscillation of the
solution. It is applied in the framework of a Particle-In-Cell method for solving some
Vlasov-Poisson equation. Since the numerical results are encouraging, several ways to
explore in the future may be usefull to improve some points in the scheme development.

We have seen that the use of $2\pi\eps$ as fast time within the first step of Algorithm
\ref{algo-time} may lead to an unstable simulation. In addition, even in a stable
simulation, the use of the particles mean period gives smaller errors than those obtained
with $2\pi\eps$. Therefore, we think it is important to find theoretically a more accurate
approximation of the fast time. 

It will be interesting to see if our numerical scheme preserves the two-scale asymptotic
limit, meaning that $\eps=0$ in the numerical scheme leads to a consistent discretization 
of the two-scale limit model. This remark is based on the fact that the ETD discretization
we have used is very close to an explicit discretization of the limit model in Theorem 1.1
of \cite{Fre06}.

Last but not least, the ETD-PIC schemes proposed in this paper need to be tested on other
systems where different types of stiff differential equations are to be solved.


\begin{thebibliography}{1}

\bibitem{BL91}  \textsc{C. K. Birdsall, A. B. Langdon}, Plasma physics via computer
simulation, Institute of Physics, Bristol (1991).

\bibitem{Boy01} \textsc{J. P. Boyd}, Chebyshev and Fourier Spectral Methods, 
Dover, New York (2001).

\bibitem{Cox-Matt02} \textsc{S. M. Cox, P. C. Matthews}, \emph{Exponential Time
Differencing for Stiff Systems}, J. Comput. Phys. \textbf{176} (2002), 430-455.

\bibitem{Fil-Son06} \textsc{F. Filbet, E. Sonnendr\"ucker}, \emph{Modeling and numerical
simulation of space charge dominated beams in the paraxial approximation}, Math. Models
Methods Appl. Sci. \textbf{16}-5 (2006), 763-791.

\bibitem{Fre06} \textsc{E. Fr\'enod}, \emph{Application of the averaging method to the
gyrokinetic plasma}, Asympt. Analysis \textbf{46} (2006), 1-28.

\bibitem{Fre-Sal-Son09} \textsc{E. Fr\'enod, F. Salvarani, E. Sonnendr\"ucker},
\emph{Long time simulation of a beam in a periodic focusing channel via a two-scale
PIC-method}, Math. Models Methods Appl. Sci. \textbf{19}-2 (2009), 175-197.

\bibitem{Hoc-Ost10} \textsc{M. Hochbruck, A. Ostermann}, \emph{Exponential integrators},
Acta Numer. \textbf{19} (2010), 209-286.

\bibitem{Kass-Tref05}\textsc{A.-K. Kassam, L. N. Trefethen}, \emph{Fourth-order
time-stepping for stiff PDEs}, SIAM J. Sci. Comput. \textbf{26}-4 (2005), 1214-1233.

\bibitem{TPLH10}\textsc{T. T\"uckmantel, A. Pukhov, J. Liljo, M. Hochbruck}, 
\emph{Three-Dimensional Relativistic Particle-in-Cell Hybrid Code Based on an Exponential
Integrator}, IEEE Trans. Plasma Sci. \textbf{38}-9 (2010), 2383-2389.

\end{thebibliography}
\end{document}